\documentclass[a4paper,12pt]{amsart}
\usepackage{amssymb}
\usepackage{ifthen}
\usepackage{graphicx}

\newtheorem{thm}{Theorem}
\newtheorem{cor}{Corollary}
\newtheorem{lem}{Lemma}

\newtheorem{rem}{Remark}

\newtheorem{conj}{Conjecture}
\theoremstyle{definition}
\theoremstyle{example}
\newtheorem{example}[equation]{Example}
\newtheorem{prob}[equation]{Problem}

\newcounter {own}
\def\theown {\thesection       .\arabic{own}}

\newenvironment{pf}[1][]{%
 \vskip 3mm
 \noindent
 \ifthenelse{\equal{#1}{}}%
  {{\slshape Proof. }}%
  {{\slshape #1.} }%
 }%
{\qed\bigskip}

\newcounter{alphabet}
\newcounter{tmp}

\makeatletter
\newcommand{\Ref}[1]{\@ifundefined{r@#1}{}{\setcounter{tmp}{\ref{#1}}\Alph{tmp}}}
\makeatother

\newcommand{\IC}{{\mathbb C}}
\newcommand{\ID}{{\mathbb D}}

\def\be{\begin{equation}}
\def\ee{\end{equation}}

\newcommand{\bee}{\begin{enumerate}}
\newcommand{\eee}{\end{enumerate}}

\newcommand{\blem}{\begin{lem}}
\newcommand{\elem}{\end{lem}}
\newcommand{\bthm}{\begin{thm}}
\newcommand{\ethm}{\end{thm}}
\newcommand{\bcor}{\begin{cor}}
\newcommand{\ecor}{\end{cor}}
\newcommand{\beg}{\begin{example}}
\newcommand{\eeg}{\end{example}}
\newcommand{\begs}{\begin{examples}}
\newcommand{\eegs}{\end{examples}}
\newcommand{\bdefe}{\begin{defin}}
\newcommand{\edefe}{\end{defin}}
\newcommand{\bprob}{\begin{prob}}
\newcommand{\eprob}{\end{prob}}
\newcommand{\bei}{\begin{itemize}}
\newcommand{\eei}{\end{itemize}}

\newcommand{\bcon}{\begin{conj}}
\newcommand{\econ}{\end{conj}}
\newcommand{\bcons}{\begin{conjs}}
\newcommand{\econs}{\end{conjs}}
\newcommand{\bprop}{\begin{propo}}
\newcommand{\eprop}{\end{propo}}
\newcommand{\br}{\begin{rem}}
\newcommand{\er}{\end{rem}}
\newcommand{\brs}{\begin{rems}}
\newcommand{\ers}{\end{rems}}
\newcommand{\bo}{\begin{obser}}
\newcommand{\eo}{\end{obser}}
\newcommand{\bos}{\begin{obsers}}
\newcommand{\eos}{\end{obsers}}
\newcommand{\bpf}{\begin{pf}}
\newcommand{\epf}{\end{pf}}
\newcommand{\ba}{\begin{array}}
\newcommand{\ea}{\end{array}}
\newcommand{\beq}{\begin{eqnarray}}
\newcommand{\beqq}{\begin{eqnarray*}}
\newcommand{\eeq}{\end{eqnarray}}
\newcommand{\eeqq}{\end{eqnarray*}}

\newcommand{\ds}{\displaystyle}

\def\cc{\setcounter{equation}{0}   
\setcounter{figure}{0}\setcounter{table}{0}}

\newcounter{minutes}
\divide\time by 60
\newcounter{hours}
\multiply\time by 60 \addtocounter{minutes}{-\time}

\begin{document}
\bibliographystyle{amsplain}
\title[On harmonic combination of univalent functions]{On harmonic combination  of univalent functions}

\thanks{
File:~\jobname .tex,
          printed: \number\day-\number\month-\number\year,
          \thehours.\ifnum\theminutes<10{0}\fi\theminutes}

\author{M. Obradovi\'{c}}
\address{M. Obradovi\'{c}, Department of Mathematics,
Faculty of Civil Engineering,
University of Belgrade, Bulevar Kralja Aleksandra 73, 11000
Belgrade, Serbia. } \email{obrad@grf.bg.ac.rs}

\author{S. Ponnusamy$^\dagger$}
\address{S. Ponnusamy, Department of Mathematics,
Indian Institute of Technology Madras, Chennai--600 036, India.}
\email{samy@iitm.ac.in}

\keywords{Coefficient inequality, partial sums, radius of univalence, analytic, univalent and starlike
functions\\
$^\dagger$ Corresponding author}
\date{\today  
;  File: opon5-10.tex}

\begin{abstract}
Let ${\mathcal S}$ be the class of all functions $f$ that are analytic and univalent in
the unit disk $\ID$ with the normalization $f(0)=f'(0)-1=0$.
Let $\mathcal{U} (\lambda)$ denote the set of all $f\in {\mathcal S}$ satisfying
the condition
$$\left |f'(z)\left (\frac{z}{f(z)} \right )^{2}-1\right | <\lambda ~\mbox{ for $z\in \ID$},
$$
for some $\lambda \in (0,1]$. In this paper, among other things, we study a ``harmonic mean'' of two univalent analytic
functions. More precisely, we discuss the properties of the class of functions $F$ of the form
$$\frac{z}{F(z)}=\frac{1}{2}\left( \frac{z}{f(z)}+\frac{z}{g(z)}
\right),
$$
where $f,g\in \mathcal{S}$ or $f,g\in \mathcal{U}(1)$. In  particular,
we determine the radius of univalency of $F$, and propose two conjectures concerning the univalency of
$F$.
\end{abstract}

\thanks{The work of the first author was supported by MNZZS Grant, No. ON174017, Serbia.
}

\maketitle
\pagestyle{myheadings}
\markboth{M. Obradovi\'{c} and S.Ponnusamy}{
Combination of univalent functions}
\cc

\section{Introduction and Main Results}\label{5-10sec1}

For each $r>0$, we denote by $\ID _r$ the open disk $\{ z\in \IC:\,|z|< r\}$ and by
$\ID$ the unit disk $\ID _1$.  Let ${\mathcal A}$ be the class of all functions $f$
that are analytic in $\ID$ with the normalization $f(0)=f'(0)-1=0$.
Denote by ${\mathcal S}$, ${\mathcal S}^*$,  ${\mathcal K}$, and ${\mathcal C}$,
the subfamilies of $\mathcal A$ that are, respectively, univalent, starlike, convex, and close-to-convex
in $\ID$ (see \cite{Du,Go} for some detailed discussion on these classes).
It is well-known that a function $f\in {\mathcal S}$ is starlike if
$$ {\rm Re}\left(\frac{zf'(z)}{f(z)} \right)>0\quad \mbox{for all $z\in \ID$}.
$$
Similarly, a function $f\in {\mathcal S}$ is close-to-convex if
$${\rm Re \,}\left ( \frac{zf'(z)}{g(z)} \right )>0 \quad \mbox{for all $z\in \ID$}
$$
for some $g\in {\mathcal S}^*$.  In \cite{Mitrino-79}, Mitrinovi\'c essentially
investigated certain geometric properties of the functions $f$ of the form
\be\label{5-10eq0}
f(z)=\frac{z}{\phi (z)}, \quad \phi (z)=1+\sum_{n=1}^{\infty}b_nz^n.
\ee
In \cite{RST-84}, Reade et al. derived coefficient conditions that guarantee the
univalence, starlikeness, or convexity of rational functions of the form
(\ref{5-10eq0}). These results have been improved and generalized
by the authors in \cite{OPSV}. In connection with a problem due to \cite{Hayman-67},
several authors (eg. \cite{Rus-Wirth76}) discussed the univalency of functions in the set of convex linear combinations
of the form
$$ \mu f(z) +(1-\mu)g(z), ~\mu \in [0,1],
$$
when $f,g$ belonging to suitable subsets of $\mathcal S$. In this paper, we shall consider a similar problem
for univalent functions $f$ of the form (\ref{5-10eq0}).

Let $\mathcal{U} (\lambda)$ denote the set of all  $f\in {\mathcal A}$ in $\ID$ satisfying
the condition (\cite{OP-01,obpo-2007a})
\be\label{4-10eq0}
\left |f'(z)\left (\frac{z}{f(z)} \right )^{2}-1\right | <\lambda,
 ~\mbox{ for $z\in \ID$}
\ee
and for some $\lambda \in (0,1]$. Functions in  $\mathcal{U} (1)=: \mathcal{U}$ is known
to be univalent in $\ID$, see \cite {Aks58,OzNu}. Clearly,  $\mathcal{U} (\lambda)\subset  \mathcal{U}$ for $\lambda \in (0,1]$
and so, functions in  $\mathcal{U} (\lambda)$ are univalent in $\ID$. Set
$$\mathcal{U}_2(\lambda)=\{f\in {\mathcal U}(\lambda):\, f''(0)=0\}.
$$
For convenience, we may also let $\mathcal{U}_2={\mathcal U}_2(0)$. It is known (eg. \cite{obpo-2007a}) that
functions in  $\mathcal{U}_2$ are included in the class ${\mathcal P}(1/2)$, where
$${\mathcal P}(1/2)=\left \{f\in {\mathcal A}:\, {\rm Re\,}(f(z)/z)>1/2 ~\mbox{ for $z\in \ID$}\right \}.
$$
We remark that ${\mathcal K} \subset {\mathcal P}(1/2)$ and there exist functions $f$ in ${\mathcal S}$
such that $f\not\in {\mathcal U}$.

It is convenient to say that  $f$ belongs to ${\mathcal U}(\lambda )$ in the disk $|z|<r$
if the inequality in (\ref{4-10eq0}) holds for $|z|< r$ instead of the whole unit disk $\ID$.
For instance, if $\lambda =1$, this is equivalent to saying that $g$ defined by
$g(z)=r^{-1}f(rz)$ belongs to $\mathcal{U}$, whenever $f$ belongs to ${\mathcal U}$ in the disk $|z|<r$.
Similar terminology will be followed for other related classes of functions, eg., starlike functions in $|z|<r$.

Now, we state our main results.

\bthm\label{5-10th1}
Let $f,g\in \mathcal{S}$. Suppose that  $\frac{f(z)+g(z)}{z}\neq 0$ for $z\in \ID$ and consider the function $F$
defined by
\be\label{5-10eq1}
F(z)=\frac{2f(z)g(z)}{f(z)+g(z)}.
\ee
Then $G$, defined by $G(z)=r^{-1}F(rz)$, belongs to $\mathcal{U}(\lambda )$ for $0<r\leq \sqrt{\lambda /(1+\lambda)}$.
In particular, $F$ belongs to $\mathcal{U}$ in the disk $|z|<1/\sqrt{2}  \approx 0.707107$ $($and hence,
$F$ is univalent in $\ID_{1/\sqrt{2}})$. In addition, $r^{-1}F(rz)$ belongs $\mathcal{S}^*$
for
$$0<r\leq r_0=\sqrt{1-(2 /(4-|b_1+c_1|))},
$$
where $b_1+c_1 =-(f''(0)+g''(0))/2=-F''(0) $.
\ethm

If $b_1+c_1=0$, then from Theorem \ref{5-10th1} we obtain that $G$ defined $G(z)=r^{-1}F(rz)$
belongs to $\mathcal{U}\cap \mathcal{S}^*$ whenever $0<r\leq 1/\sqrt{2}$. Moreover,
since $\mathcal{U}\subsetneq \mathcal{S}$, it is natural to prove an analog of Theorem \ref{5-10th1} by replacing
the assumption $f,g\in\mathcal{S} $ by $f,g\in\mathcal{U}.$ Now, we are in a position to state our next result.

\bthm\label{5-10th3}
Let $f\in \mathcal{U}(\lambda _1 )$, $g\in \mathcal{U}(\lambda _2 )$ $(0<\lambda _1 ,\, \lambda _2\leq 1)$
and  $\frac{f(z)+g(z)}{z}\neq 0$ for $z\in \ID$.
Define $F$ by {\rm (\ref{5-10eq1})} and $G$  by $G(z)=r^{-1}F(rz)$.
Then $G$ belongs to $\mathcal{U}(\lambda )$ whenever
\be\label{5-10eq6}
0<r\leq \sqrt{\frac{-K^2+K\sqrt{K^2+4}}{2}}
~\mbox{ with }~ K=\sqrt{2\lambda ^2/(\lambda _1 +\lambda _2)}.
\ee
In particular, if $f,g\in \mathcal{U}$, then $G\in \mathcal{U}$ for $0<r\leq \sqrt{\frac{\sqrt{5}-1}{2}}$; that is
$F$ is univalent in the disk $|z|<\sqrt\frac{\sqrt{5}-1}{2}\approx 0.78615$.
\ethm

At this place, it is appropriate to present a two parameters family of analytic functions dealing with a number
of issues concerning our investigation.

\vspace{8pt}
\noindent
{\bf Example 1.~}  For $0\neq \alpha \in [-1,1]$, we consider
$$f_\alpha (z)= \frac{z(1-\alpha z)}{1-z^{2}}.
$$
By a computation, we obtain that
$$(1-z^{2})f_\alpha ' (z)= \frac{1+z^{2}-2\alpha z}{1-z^{2}}
$$
and so,
$$ {\rm Re\,}((1-z^{2})f_\alpha ' (z))= \frac{(1-|z|^2)(1+|z|^2-2\alpha{\rm Re\,} z)}{|1-z^{2}|^2} >0 ~\mbox{ for $z\in \ID$}.
$$
We conclude that for each $\alpha$, the function $f_\alpha$ is close-to-convex in $\ID$.
Now, let $F(z)=F_{\alpha ,\beta}(z)$ in the unit disk $\ID$ be defined by
$$\frac{z}{F (z)}=\frac{1}{2}\left( \frac{z}{f_\alpha(z)}+\frac{z}{f_{\beta}(z)}\right)
$$
where $\alpha, \beta \in [-1,1]\backslash \{0\}$. A computation gives
$$F(z)=\frac{z(1-\alpha z)(1-\beta z)}{(1-z^2)(1-((\alpha +\beta)/2)z)}
$$
and
$$\frac{z}{F(z)}=1+\sum_{n=1}^{\infty}b_nz^n,
$$
where
$$b_1=\frac{\alpha +\beta}{2} ~\mbox{ and }~b_n=-\frac{1}{2}\left ( \alpha ^{n-2}(1-\alpha ^2)+\beta ^{n-2}(1-\beta ^2)
\right ) ~\mbox{ for $n\geq 2$}.
$$
First we wish to show that
$$S:=\sum_{n=2}^{\infty}(n-1)|b_n|^2\leq 1
$$
which is a necessary condition for $F$ to belong to ${\mathcal S}$ (see the well-known Area
Theorem \cite[Theorem 11 on p.193 of Vol. 2]{Go}). As
$$4|b_n|^2= \alpha ^{2(n-2)}(1-\alpha ^2)^2+\beta ^{2(n-2)}(1-\beta ^2)^2+2 (1-\alpha ^2)(1-\beta ^2)(\alpha \beta )^{n-2}
$$
for $n\ge 2$ and
$$\sum_{n=2}^{\infty}(n-1)x^{n-2}= \frac{1}{(1-x)^2} ~\mbox{ for $|x|<1$},
$$
it follows easily that
$$S=\frac{1}{4}\left (1+1+\frac{2 (1-\alpha ^2)(1-\beta ^2)}{(1-\alpha \beta)^2}\right )\leq 1
$$
and the equality holds if $\alpha =\beta$. Thus, $F$ satisfies the necessary condition for $F$
to belong to the class ${\mathcal S}$ whenever  $\alpha, \beta \in [-1,1]\backslash \{0\}$.
On the other hand for certain values of $\alpha, \beta $ these functions $F=F_{\alpha ,\beta}$
belong neither to $\mathcal U$ nor to ${\mathcal S}^*$.

We note that $z/F(z) \neq 0$ in $\ID$ and, by a lengthy computation, we obtain that
$$\left(\frac{z}{F(z)}\right)^{2}F'(z)-1 =z^2\frac{(1-\alpha\beta)(1-\alpha z)(1-\beta z) -
(1-z^2)((\alpha-\beta )^2/2)}{(1-\alpha z)^2(1-\beta z)^2}
$$
Setting $\beta =-\alpha,$ we see that for the function  $F_\alpha(z):=F_{\alpha ,-\alpha}(z)$, we have
$$\left(\frac{z}{F_\alpha (z)}\right)^{2}F_\alpha '(z)-1 =-(1-\alpha ^2)\frac{z^2(1+\alpha ^2z^{2})}{(1-\alpha ^2z^2)^2}.
$$
Therefore, we have
$$\left|\left(\frac{z}{F_\alpha (z)}\right)^{2}F_\alpha '(z)-1\right|
\leq  (1-\alpha ^2 )\frac{|z|^2(1+\alpha ^2|z|^2)}{(1-\alpha ^2|z|^2)^2}
$$
which is less than $1$ whenever
$$(2\alpha ^2-1 )\alpha ^2|z|^4 - (1+\alpha ^2)|z|^2 +1 >0.
$$
Solving the last inequality gives the condition $|z|^2 <r^2_{\mathcal U}$, where
$$ r_{\mathcal U}=\sqrt{\frac{2}{1+\alpha ^2+\sqrt{(1-\alpha ^2)(7\alpha ^2 +1)} }}.
$$
The above discussion shows that
$$\left|\left(\frac{z}{F_\alpha (z)}\right)^{2}F_\alpha '(z)-1\right| <1 ~\mbox{ for $|z|< r_{\mathcal U}$.}
$$
Also, for $z =r$, where $r_{\mathcal U}\leq r<1$, we have that
$$ \left|\left(\frac{z}{F_\alpha (z)}\right)^{2}F_\alpha '(z)-1\right|=
(1-\alpha ^2 )\frac{r^2(1+\alpha ^2r^2)}{(1-\alpha ^2r^2)^2} \geq 1 ,
$$
showing that the function $F_\alpha $ is in the class $\mathcal{U}$ in the disk
$|z|<r_{\mathcal U}$ (so $F_\alpha $  is univalent in this disk) but not in any larger disk.
That is, $r^{-1}F_\alpha (rz)$ belongs to $\mathcal{U}$ for $0<r\leq r_{\mathcal U}$, but not
for a larger value of $r$.

Next, we show that for certain values of $\alpha, \beta $, the functions
$F=F_{\alpha, \beta }$ are not starlike in the unit disk $\ID$. A straightforward
computation shows that
$$F'(z)=\frac{M(z)}{(1-z^2)^2(1-((\alpha +\beta)/2)z)^2},
$$
where
\beqq
M(z)& =& 1-2(\alpha +\beta) z +(1+3\alpha\beta +(\alpha +\beta)^2/2)z^2\\
& & ~~~~~~ -(\alpha +\beta)(1+\alpha\beta)z^3 +((\alpha ^2 +\beta^2)/2)z^4.
\eeqq
If $\beta=-\alpha $ with $|\alpha |>1/9$ then in this case $M(z)$ takes the form
$$M(z)=\alpha ^2[z^2+A][z^2+\overline{A}], \quad A=\frac{1-3\alpha ^2+i\sqrt{(9\alpha ^2-1)(1-\alpha^2)}}{2\alpha ^2}
$$
and we see that $|A|\geq 1$ showing that $F_\alpha '(z)\neq 0$ in $\ID$.

Also, if $|\alpha |\leq 1/9$, then we see that
$$M(z)=\alpha ^2[z^2+B_{+}][z^2+ B_{-}]$$
where
$$B_{\pm}=\frac{1-3\alpha ^2\pm\sqrt{(1-9\alpha ^2)(1-\alpha^2)}}{2\alpha ^2}\geq 1 .
$$
Again, we see that $F_\alpha'(z)\neq 0$ in $\ID$. Thus,  $F_\alpha$ is locally
univalent in $\ID$.

On the other hand,  it follows easily that
$$ \frac{zF_\alpha'(z)}{F_\alpha(z)}=
\frac{1+(1-3\alpha ^2)z^2+\alpha ^2z^4}{(1-\alpha ^2z^2)(1-z^2)}.
$$
A straightforward computation shows that for $0<\theta <\pi$, 
$$ {\rm Re\,}\left (\frac{e^{i\theta}F_\alpha'(e^{i\theta})}{F_\alpha(e^{i\theta})}\right )=
\frac{A(\theta )}{|1-\alpha ^2e^{2i\theta}|^2\,|1-e^{2i\theta}|^2},
$$
where
$$A(\theta )=4\alpha ^2(\alpha ^2-\cos 2\theta) (1-\cos 2\theta).
$$
Therefore, $A(\theta )<0$
if $0<\alpha ^2 <\cos 2\theta <1$, i.e. $|\theta|<(1/2)\arccos (\alpha ^2) <\pi /4$.
This observation shows that for each $\alpha \in (-1,1)\backslash \{0\}$, the function $F_\alpha $
is not starlike in $\ID$ although $F_\alpha $ is locally univalent in $\ID$.


\vspace{8pt}

The above example motivates the following conjectures

\bcon
\begin{enumerate}
\item[{\rm \textbf{(a)}}]
The function $F$ defined by {\rm (\ref{5-10eq1})} is not necessarily univalent in $\ID$ whenever
$f,g\in \mathcal{S}$ such that  $((f(z)+g(z))/z) \neq 0$ in  $\ID$.
\item[{\rm \textbf{(b)}}] The function $F$ defined by {\rm (\ref{5-10eq1})} is  univalent in $\ID$ whenever
$f,g\in \mathcal{C}$ such that  $((f(z)+g(z))/z) \neq 0$ in $\ID$.
\end{enumerate}
\econ

Theorem \ref{5-10th1} may be generalized in the following form.

\bthm\label{5-10th2}
Let $f_k\in \mathcal{S}$ for $k=1,\ldots, m$ and $\sum_{k=1}^{m}\frac{z}{f_{k}(z)}\neq 0$ for $z\in\ID$.
Define $F$ by
\be\label{5-10eq4}
\frac{z}{F(z)}=\frac{1}{m} \sum_{k=1}^{m}\frac{z}{f_{k}(z)}.
\ee
Then we have
\begin{enumerate}
\item[\textbf{(a)}] $G$ defined by $G(z)=r^{-1}F(rz)$ belongs to $\mathcal{U}(\lambda )$
for $0<r\leq \sqrt{\lambda /(1+\lambda)}$. In particular, $F$ is univalent in the disk $|z|<1/\sqrt{2}$.
\item[\textbf{(b)}]  $G$ belongs to $\mathcal{S}^*$
for $0<r\leq \sqrt{\lambda /(1+\lambda)}$, with
$$\lambda =1- \frac{1}{m} \left |\sum_{k=1}^{m} \frac{f_k''(0)}{2} \right |.
$$
In particular, $F$ is starlike $($univalent$)$ in the disk $|z|<1/\sqrt{2}$ whenever
$f_k''(0)=0$ for each $k=1,\ldots , m$.
\end{enumerate}
\ethm

The idea of the proof of Theorem \ref{5-10th3} can be used to prove the following general result. We omit its
proof.

\bthm\label{5-10th4}
Let $f_k\in \mathcal{U}(\lambda _k )$ $(0<\lambda _k\leq 1)$ for $k=1,\ldots, m$,
$\sum_{k=1}^{m}\frac{z}{f_{k}(z)}\neq 0$ for $z\in\ID$
and $F$ be defined by {\rm (\ref{5-10eq4})}.
Then  $G$ defined by $G(z)=r^{-1}F(rz)$ belongs to $\mathcal{U}(\lambda )$ whenever
\be\label{5-10eq7}
0<r\leq \sqrt{\frac{-K^2+K\sqrt{K^2+4}}{2}}
~\mbox{ with }~ K=\sqrt{\frac{m\lambda ^2}{\sum_{k=1}^{m}\lambda _k}}.
\ee
In particular, if $f_k\in \mathcal{U}$ for $k=1,\ldots, m$,
then $G\in \mathcal{U}$ for $0<r\leq \sqrt{\frac{\sqrt{5}-1}{2}}$; that is
$F$ is univalent in the disk $|z|<\sqrt\frac{\sqrt{5}-1}{2}$.
\ethm

The proof of Theorems \ref{5-10th1}, \ref{5-10th3}, \ref{5-10th2}, and  \ref{5-10th4} are presented in Section \ref{5-10sec3}.

\section{Preliminary Lemmas}\label{5-10sec2}

For the proofs of our results, we need the following lemmas.

\blem\label{5-10lem1}
Let $\phi(z)=1+\sum_{n=1}^\infty b_nz^n$ be a non-vanishing analytic
function on $\ID$ and let $f$ be of the form {\rm (\ref{5-10eq0})}.
Then, we have the following:
\begin{enumerate}
\item[\textbf{(a)}] If
$\sum_{n=2}^\infty (n-1)|b_n|\leq \lambda ,$ then $f\in {\mathcal U}(\lambda )$.
\item[\textbf{(b)}] If $\sum_{n=2}^\infty (n-1)|b_n|\leq 1-|b_1|$, then $f\in {\mathcal S}^*$.
\item[\textbf{(c)}] If $f\in {\mathcal U}(\lambda )$, then $\sum_{n=2}^{\infty}(n-1)^2|b_n|^2\leq \lambda ^2.$
\end{enumerate}
\elem

The conclusion \textbf{(a)} in Lemma \ref{5-10lem1} is from \cite{OP-01,OPSV} whereas the
\textbf{(b)} is due to Reade et al. \cite[Theorem 1]{RST-84}. Finally, as  $f\in {\mathcal U}(\lambda)$, we have
$$\left |f'(z)\left (\frac{z}{f(z)} \right )^{2}-1\right |=
\left |-z \left ( \frac{z}{f(z)} \right )'+\frac{z}{f(z)} -1\right | = \left |\sum_{n=2}^{\infty}(n-1)b_{n}z^{n}\right |
\leq \lambda
$$
and so \textbf{(c)} follows from Prawitz' theorem which is an immediate consequence of Gronwall's area theorem.
This may be also obtained as a consequence of Parseval's  relation.

Next we recall the following result due to Obradovi\'c and Ponnusamy \cite{obpo-2009a}.

\blem\label{5-10lem2}
Let $f\in {\mathcal A}$ have the form
\be\label{5-10eq8}
\frac{z}{f(z)}=1+b_{1}z+b_{2}z^{2}+ \cdots ~ \mbox{ with $b_{n}\geq 0$ for all $n\ge 2$}
\ee
and for all $z$ in a neighborhood of $z=0$. Then the following conditions are equivalent.
\bee
\item[\textbf{(a)}] $f\in\mathcal{S}$,
\item[\textbf{(b)}] $\ds \frac{f(z)f'(z)}{z} \neq 0$ for $z\in \ID$,
\item[\textbf{(c)}] $\ds \sum_{n=2}^{\infty}(n-1)b_{n}\leq 1$,
\item[\textbf{(d)}] $f\in\mathcal{U}$.
\eee
\elem

This lemma helps to compare the relation between results here and the earlier
work of the authors in \cite{obpo-2009a}, in particular.

\section{Proofs}\label{5-10sec3}

\noindent \textbf{Proof of Theorem \ref{5-10th1}.}
Let $f,g\in \mathcal{S}$. Then $f$ and $g$ can be written in the form
\be\label{5-10eq5}
\frac{z}{f(z)}=1+b_{1}z+b_{2}z^{2}+ \cdots   ~\mbox{ and }~ \frac{z}{g(z)}=1+c_{1}z+c_{2}z^{2}+ \cdots .
\ee
Further, as $f,g \in \mathcal{S}$, the well-known Gronwall's Area Theorem \cite[Theorem 11 on p.193 of Vol. 2]{Go}
gives
\be\label{5-10eq2}
\sum_{n=2}^\infty (n-1)|b_n|^2\leq 1  ~\mbox{ and }~\sum_{n=2}^{\infty}(n-1)|c_{n}|^{2}\leq1.
\ee
From (\ref{5-10eq1}), we may rewrite $F$ in the form
$$\frac{z}{F(z)}=\frac{1}{2}\left( \frac{z}{f(z)}+\frac{z}{g(z)}
\right)=1+\sum_{n=1}^{\infty}\frac{b_{n}+c_{n}}{2}z^{n}.
$$
For $0<r\leq 1$, we define $G$ by $G(z)=r^{-1}F(rz)$ so that
$$\frac{z}{G(z)}=\frac{z}{r^{-1}F(rz)}=1+\sum_{n=1}^{\infty}\frac{b_{n}+c_{n}}{2}r^{n}z^{n}.
$$
In order to prove that $F$ is univalent in $|z|<1/\sqrt{2}$, it suffices to show that
$G\in\mathcal U$ for $0<r\leq 1/\sqrt{2}$. According to Lemma \ref{5-10lem1}\textbf{(a)}
(compare with Lemma \ref{5-10lem2}\textbf{(c)}),
it suffices to show that
\be\label{5-10eq3}
S:=\sum_{n=2}^{\infty}(n-1)\left|\frac{b_{n}+c_{n}}{2}\right|r^{n}\leq 1
\ee
for $0<r\leq 1/\sqrt{2}$. By (\ref{5-10eq2}) and the Cauchy-Schwarz inequality, we have
$$\sum_{n=2}^{\infty}(n-1)|b_{n}|r^{n}
\leq \left(\sum_{n=2}^{\infty}(n-1)|b_{n}|^{2}\right)^{\frac{1}{2}}\left(\sum_{n=2}^{\infty}(n-1)r^{2n}\right)^{\frac{1}{2}}
\leq \frac{r^{2}}{1-r^{2}}.
$$
and similarly, we obtain that
$$\sum_{n=2}^{\infty}(n-1)|c_{n}|r^{n}\leq  \frac{r^{2}}{1-r^{2}}.
$$
As $|b_{n}+c_{n}|\leq |b_{n}|+|c_{n}|$, the last two inequalities gives that
$$S\leq \frac{r^{2}}{1-r^{2}}.
$$
Thus, $S\leq 1$ whenever $r^{2}\leq 1-r^{2} $, i.e. if $r\leq 1/\sqrt{2}.$ Thus, $G\in {\mathcal U}$
and we complete the proof of the first part. The proof of the second part is a consequence of
Lemma \ref{5-10lem1}\textbf{(a)} and solving the inequality $r^{2}\leq \lambda (1-r^{2})$. The final part,
namely, $G\in \mathcal{S}^*$, follows by setting $\lambda =1-|b_1+c_1|/2$ and applying Lemma \ref{5-10lem1}\textbf{(b)}.
In other words, $F(|z|<r_0)$ is a starlike domain, where
$r_0=\sqrt{\lambda /(1+\lambda)}$ with $\lambda =1-|b_1+c_1|/2$. A computation gives

\vspace{6pt}
\hfill
$r_0=\sqrt{1-(2 /(4-|b_1+c_1|))}.
$
\hfill $\Box$

\vspace{9pt}

\noindent \textbf{Proof of Theorem \ref{5-10th3}.}
Let $f\in \mathcal{U}(\lambda _1 )$ and $g\in \mathcal{U}(\lambda _2 )$, and have the form
(\ref {5-10eq5}). By Lemma \ref{5-10lem1}\textbf{(c)}, we have
\be\label{5-10eq9a}
\sum_{n=2}^{\infty}(n-1)^{2}|b_{n}|^{2}\leq \lambda_1^2 ~\mbox{ and }
\sum_{n=2}^{\infty}(n-1)^{2}|c_{n}|^{2}\leq \lambda_2^2.
\ee
As in the proof of Theorem \ref{5-10th1}, for $G$ belonging to  $\mathcal{U}(\lambda )$,
it suffices to show by Lemma \ref{5-10lem1}\textbf{(a)} that
$$T=\sum_{n=2}^{\infty}(n-1)\left|\frac{b_{n}+c_{n}}{2}\right|r^{n}\leq \lambda
$$
under the condition (\ref{5-10eq6}).  Now, by \eqref{5-10eq9a} and the Cauchy-Schwarz inequality, we have
$$\sum_{n=2}^{\infty}(n-1)|b_{n}|r^{n}
\leq \left(\sum_{n=2}^{\infty}(n-1)^2|b_{n}|^{2}\right)^{\frac{1}{2}}\left(\sum_{n=2}^{\infty}r^{2n}\right)^{\frac{1}{2}}
\leq \frac{\lambda _1 r^{2}}{\sqrt{1-r^{2}}},
$$
and similarly, we obtain that
$$\sum_{n=2}^{\infty}(n-1)|c_{n}|r^{n}\leq  \frac{\lambda _2 r^{2}}{\sqrt{1-r^{2}}}.
$$
In view of the last two inequalities, it follows that
$$T\leq  \left (\frac{\lambda _1 +\lambda _2}{2}\right )\frac{r^{2}}{\sqrt{1-r^{2}}}.
$$
It can be easily seen that the last expression is less than or equal to $\lambda$ if
and only if $r$ satisfies the  inequality (\ref{5-10eq6}).
This means that the function $G\in \mathcal{U}(\lambda )$ under the condition (\ref{5-10eq6}),
which is equivalent to saying that $F\in \mathcal{U}$ (and hence univalent) in the disk
$$|z|< \sqrt{\frac{-K^2+K\sqrt{K^2+4}}{2}}, \quad K=\sqrt{\frac{2\lambda ^2}{\lambda _1 +\lambda _2}}.
$$
The proof of the main part is complete. Setting $\lambda_1=\lambda_2=\lambda =1$, it follows that
if $f,g\in \mathcal{U}$, then $G\in \mathcal{U}$ for $0<r\leq \sqrt{\frac{\sqrt{5}-1}{2}}$. In particular,
we obtain that $F$ is univalent in the disk $|z|<\sqrt\frac{\sqrt{5}-1}{2}\approx 0.78615$.
\hfill $\Box$

\vspace{9pt}

\noindent \textbf{Proof of Theorem \ref{5-10th2}.}
As $f_k\in \mathcal{S}$ for $k=1,\ldots, m$, we may represent $z/f_k(z)$ in power series form
$$\frac{z}{f_k(z)}=1+\sum_{n=1}^{\infty} b_{n}^{(k)}z^n
$$
so that
$$\frac{z}{F(z)}=1+\sum_{n=1}^{\infty} B_{n}^{(m)}z^n, \quad B_n^{(m)}=\frac{1}{m} \sum_{k=1}^{m}b_{n}^{(k)}.
$$

As the remaining part of the proof follows in the same lines as that of the proof of
Theorem \ref{5-10th1}, we omit its details. For example, for part \textbf{(b)},
the role of $\lambda =1-|b_1+c_1|/2$ in the proof of Theorem \ref{5-10th1}
will be replaced by
$$\lambda =1- \frac{1}{m} \left |\sum_{k=1}^{m} \frac{f_k''(0)}{2} \right |.
$$
\hfill $\Box$

\section{Discussion}\label{5-10sec4}

In Theorem \ref{5-10th4},  it is possible to remove the hypothesis that
$\sum_{k=1}^{m}\frac{z}{f_{k}(z)}\neq 0$ for $z\in\ID$. For example, if $f_k\in \mathcal{U}(\lambda _k  )$
 with $f_k''(0)=0$ for $k=1,\ldots, m$, then it is known that \cite{OPSV,obpo-2007a}
$${\rm Re\,} \left (\frac{f_k(z)}{z}\right ) >\frac{1}{1+\lambda _k}
\geq \frac{1}{2}\quad \mbox{for $z\in \ID$} .
$$
In particular, we obtain that ${\rm Re\,} (f_k(z)/z)>0$ in $\ID$ and for each $k=1,\ldots, m$.
Thus, ${\rm Re\,} (z/f_k(z))>0$ in $\ID$ and for each $k=1,\ldots, m$ so that the assumption that
$\sum_{k=1}^{m}\frac{z}{f_{k}(z)}\neq 0$ obviously holds for $z\in\ID$.
In this case, this observation gives that $F$ defined by (\ref{5-10eq4}) is univalent and starlike in the disk
$|z|<\sqrt\frac{\sqrt{5}-1}{2}$ whenever $f_k\in \mathcal{U}$ with $f_k''(0)=0$ for $k=1,\ldots, m$.
Moreover, in some special situations, one can improve Theorem \ref{5-10th1}. For instance, we have

\bthm\label{5-10th5}
Let $f,g\in {\mathcal S}$ have the form {\rm (\ref{5-10eq8})} with $b_{n}\geq 0$  and $c_{n}\geq 0$ for all $n\ge 2$.
Then the function $F$ defined by {\rm (\ref{5-10eq1})} belongs to ${\mathcal S}$. In particular, if $b_1+c_1=0$, then $F$ is
also starlike in $\ID$.
\ethm
\bpf
By Lemma \ref{5-10lem2}, we have $f,g\in\mathcal{U}$ and
$$\sum_{n=2}^{\infty}(n-1)b_{n}\leq 1 ~\mbox{ and }~
\sum_{n=2}^{\infty}(n-1)c_{n}\leq 1.
$$
The last two coefficient conditions imply that
$$\sum_{n=2}^{\infty}(n-1)\frac{b_{n}+c_{n}}{2} \leq 1
$$
and therefore, $F$ defined by {\rm (\ref{5-10eq1})} is univalent in  $\ID$.

If $b_1+c_1=0$, then according to Lemma \ref{5-10lem1}\textbf{(b)} the function $F$
defined by {\rm (\ref{5-10eq1})} is starlike in $\ID$.
\epf

Using Theorem \ref{5-10th5}, one may state a general result as in Theorems \ref{5-10th2} and \ref{5-10th4}.
Also, it is an open problem to determine the exact radii in Theorems \ref{5-10th1} and \ref{5-10th3}.

\end{document}